\documentclass[a4paper,11pt]{amsart}
\usepackage{amssymb,amscd,verbatim, amsthm,amsmath,amsgen,xspace}
\usepackage{latexsym}
\usepackage{amsfonts}
\usepackage{color}
\usepackage[all]{xy}

\newcommand \ch[1]{{\check{#1}}}
\newcommand \bb[1]{{\mathbb #1}}

\newcommand \R{{\bb R}}
\newcommand \Z{{\bb Z}}

\newcommand\Hom{\operatorname{Hom}}

\def\ch{\mathop{\hbox {ch}}\nolimits}

\def\End{\mathop{\hbox {End}}\nolimits}
\def\Hom{\mathop{\hbox {Hom}}\nolimits}

\DeclareMathOperator\im{Im}

\def\ker{\mathop{\hbox{Ker}}\nolimits}

\newcommand{\del}{\partial}

\newcommand{\ra}{\rightarrow}
\newcommand{\lra}{\longrightarrow}

\newcommand{\pf}{\begin{proof}}
\newcommand{\epf}{\end{proof}}
\newcommand{\eq}{\begin{equation}}
\newcommand{\eeq}{\end{equation}}
\newcommand{\eqn}{\begin{equation*}}
\newcommand{\eeqn}{\end{equation*}}

\newcommand{\frb}{\mathfrak{b}}

\newcommand{\frg}{\mathfrak{g}}
\newcommand{\frh}{\mathfrak{h}}

\newcommand{\frn}{\mathfrak{n}}
\newcommand{\fro}{\mathfrak{o}}
\newcommand{\frp}{\mathfrak{p}}

\newcommand{\frs}{\mathfrak{s}}

\newcommand{\frsp}{\mathfrak{sp}}

\newcommand{\bbC}{\mathbb{C}}

\newcommand{\bbR}{\mathbb{R}}

\newcommand{\sgn}{\operatorname{sgn}}

\newtheorem{thm}[equation]{Theorem}

\newtheorem{defi}[equation]{Definition}

\numberwithin{equation}{section}

\let\ssize\scriptstyle
\newif\ifFIRST\newdimen\MAXright\MAXright0pt
\def\sdynkin{\bgroup\eightpoint\dynkin}
\def\endsdynkin{\enddynkin\egroup}
\def\dynkin{\bgroup\FIRSTtrue\hskip.5em\setbox1\hbox{$\diagup$}%
	\setbox2\hbox{$\diagdown$}%
	\setbox0\hbox to2\wd1{\hrulefill}%
	\setbox3\hbox{$\bullet$}%
	\setbox4\hbox{$\times$}%
	\setbox7\hbox{$\circ$}%       (L.K.)
	\def\whiteroot##1{\ifFIRST\setbox5\hbox{$##1$}\ifdim\wd5>1.3em%       (L.K.)
		\hskip-.5em\hskip.5\wd5\fi\fi\FIRSTfalse%                             (L.K.)
		\hskip-.25em\raise1.5\wd3\hbox to0pt{\hss\hskip.45em$%                (L.K.)
			\ssize##1$\hss}\copy7\hskip-.25em\setbox6\hbox{$##1$}%                (L.K.)
		\MAXright\wd6}%                                                       (L.K.)
	\def\root##1{\ifFIRST\setbox5\hbox{$##1$}\ifdim\wd5>1.3em%
		\hskip-.5em\hskip.5\wd5\fi\fi\FIRSTfalse%
		\hskip-.25em\raise1.5\wd3\hbox to0pt{\hss\hskip.45em$%
			\ssize##1$\hss}\copy3\hskip-.25em\setbox6\hbox{$##1$}%
		\MAXright\wd6}%
	\def\whitedroot##1{\ifFIRST\setbox5\hbox{$##1$}\ifdim\wd5>1.3em% (L.K.)
		\hskip-.5em\hskip.5\wd5\fi\fi\FIRSTfalse% (L.K.)
		\hskip-.25em\lower1.8\wd3\hbox to0pt{\hss\hskip.45em$%  (L.K.)
			\ssize##1$\hss}\copy7\hskip-.25em\setbox6\hbox{$##1$}% (L.K.)
		\MAXright\wd6}%
	\def\whiterroot##1{\hskip-.25em\copy7\hbox to0pt{\hskip.3em$\ssize##1$\hss}%
		\hskip-.25em\setbox6\hbox{\hskip.6em$##1##1$}%
		\MAXright\wd6}%
	\def\droot##1{\ifFIRST\setbox5\hbox{$##1$}\ifdim\wd5>1.3em%
		\hskip-.5em\hskip.5\wd5\fi\fi\FIRSTfalse%
		\hskip-.25em\lower1.8\wd3\hbox to0pt{\hss\hskip.45em$%
			\ssize##1$\hss}\copy3\hskip-.25em\setbox6\hbox{$##1$}%
		\MAXright\wd6}%
	\def\rroot##1{\hskip-.25em\copy3\hbox to0pt{\hskip.3em$\ssize##1$\hss}%
		\hskip-.25em\setbox6\hbox{\hskip.6em$##1##1$}%
		\MAXright\wd6}%
	\def\norroot##1{\hskip-.36em\copy4\hbox to0pt{\hskip.3em$\ssize##1$\hss}%
		\hskip-.48em\setbox6\hbox{\hskip.6em$##1##1$}%
		\MAXright\wd6}%
	\def\noroot##1{\ifFIRST\setbox5\hbox{$##1$}\ifdim\wd5>1.3em%
		\hskip-.5em\hskip.5\wd5\fi\fi\FIRSTfalse%
		\hskip-.36em\raise1.5\wd3\hbox to0pt{\hss\hskip.6em$%
			\ssize##1$\hss}\copy4\hskip-.38em\setbox6\hbox{$##1$}%
		\MAXright\wd6}%
	\def\nodroot##1{\ifFIRST\setbox5\hbox{$##1$}\ifdim\wd5>1.3em%
		\hskip-.5em\hskip.5\wd5\fi\fi\FIRSTfalse%
		\hskip-.36em\lower1.8\wd3\hbox to0pt{\hss\hskip.6em$%
			\ssize##1$\hss}\copy4\hskip-.38em\setbox6\hbox{$##1$}%
		\MAXright\wd6}%
	\def\nolink{\hskip\wd0}%      (L.K.)
	\def\link{\raise.22em\copy0}%
	\def\llink##1{\raise.32em\copy0\hskip-\wd0%
		\raise.12em\copy0\hskip-.5\wd0\hbox to0pt{\hss$##1$\hss}\hskip.5\wd0}%
	\def\lllink##1{\raise.22em\copy0\hskip-\wd0\raise.32em\copy0\hskip-\wd0%
		\raise.12em\copy0\hskip-.5\wd0\hbox to0pt{\hss$##1$\hss}\hskip.5\wd0}%
	\def\rootupright##1{\hbox to0pt{\raise.45em\copy1\hskip-.25em\raise1.3\ht1%
			\hbox{\copy3\hskip.3em$\ssize##1$}\hss}%
		\setbox6\hbox{\hskip.6em\copy1\copy1$##1##1$}%
		\ifdim\MAXright<\wd6\MAXright\wd6\fi}%
	\def\whiterootupright##1{\hbox to0pt{\raise.45em\copy1\hskip-.25em\raise1.3\ht1% (L.K.)
			\hbox{\copy7\hskip.3em$\ssize##1$}\hss}% (L.K.)
		\setbox6\hbox{\hskip.6em\copy1\copy1$##1##1$}% (L.K.)
		\ifdim\MAXright<\wd6\MAXright\wd6\fi}% (L.K.)
	\def\norootupright##1{\hbox to0pt{\raise.45em\copy1\hskip-.36em\raise1.3\ht1%
			\hbox{\copy4\hskip.3em$\ssize##1$}\hss}%
		\setbox6\hbox{\hskip.6em\copy1\copy1$##1##1$}%
		\ifdim\MAXright<\wd6\MAXright\wd6\fi}%
	\def\rootdownright##1{\hbox to0pt{\raise-.5em\copy2\hskip-.25em\raise-1.35\ht1%
			\hbox{\copy3\hskip.3em$\ssize##1$}\hss}\setbox6%
		\hbox{\hskip.6em\copy2\copy2$##1##1$}%
		\ifdim\MAXright<\wd6\MAXright\wd6\fi}%
	\def\whiterootdownright##1{\hbox to0pt{\raise-.5em\copy2\hskip-.25em\raise-1.35\ht1% (L.K.)
			\hbox{\copy7\hskip.3em$\ssize##1$}\hss}\setbox6% (L.K.)
		\hbox{\hskip.6em\copy2\copy2$##1##1$}% (L.K.)
		\ifdim\MAXright<\wd6\MAXright\wd6\fi}% (L.K.)
	\def\rootdown##1{\hbox to0pt{\hskip-.05em\vrule height.25em depth.65em%
			\hskip-.25em\raise-.95em\hbox{\copy3\hskip.3em$\ssize##1$}\hss}%
		\setbox6\hbox{$##1$}%
		\ifdim\MAXright<\wd6\MAXright\wd6\fi}%
	\def\whiterootdown##1{\hbox to0pt{\hskip-.05em\vrule height.25em depth.65em% (L.K.)
			\hskip-.25em\raise-.95em\hbox{\copy7\hskip.3em$\ssize##1$}\hss}% (L.K.)
		\setbox6\hbox{$##1$}% (L.K.)
		\ifdim\MAXright<\wd6\MAXright\wd6\fi}% (L.K.)
	\def\dots{\hskip.5em\cdots\hskip.5em}}%
\def\enddynkin{\ifdim\MAXright>1em\hskip.5\MAXright\else\hskip.5em\fi\egroup}%

\begin{document}
%\today

\title[Symplectic Dirac cohomology and lifting of characters]{Symplectic Dirac cohomology and lifting of characters to metaplectic groups}

\author{Jing-Song Huang}

\address[Huang]{Department of Mathematics, Hong Kong University of Science and Technology,
Clear Water Bay, Kowloon, Hong Kong SAR, China}
\email{mahuang@ust.hk}

\thanks{The research described in this paper is supported by grants No. 16303218  from
Research Grant Council of HKSAR}
	
\begin{abstract} 
We formulate the  transfer factor of character lifting from orthogonal groups to symplectic groups by Adams
in the framework of symplectic Dirac cohomology for the Lie superalgebras and the Rittenberg-Scheunert correspondence of
representations of the Lie superalgebra $\fro\frsp(1|2n)$ and the Lie algebra $\fro(2n+1)$.  This leads to formulation of
a direct lifting of characters  from the linear
symplectic group $Sp(2n,\bbR)$ to its nonlinear covering metaplectic group $Mp(2n,\bbR)$.
\end{abstract}

%\dedicatory{Dedicated to the memory of Kostant}
\maketitle

%%%%%%%%%%%%%%%%%%%%%%%%%%%%%%%%%%%%%%%%%%%%%%%%%%%%%%%%%%%%%%%%%%%%%%%%%%%%%%%%
\section{Introduction}\label{section intro}

A lifting of characters on orthogonal groups and nonlinear metaplectic groups over real numbers was defined and
studied by Jeff Adams \cite{A},  and it was extended to the p-adic case by Tatiana Howard \cite{Hd}. 
 This lifting is closely
 related to the endoscopic transfer in Langlands functoriality and the theta correspondence in Howe's reductive dual pairs
 and also appeared in the work of David Renard \cite{R} and Wen-Wei Li \cite{Li}.

Dirac operators were used for geometric construction of discrete series representations by Parthasarathy \cite{P}, Atiyah and Schmid \cite{AS}, and tempered representations by Wolf \cite{W}. In the late 1990's, Vogan \cite{V} made a conjecture on the algebraic property of the Dirac operators in the Lie algebra setting.  This conjecture was proved by Pandzic and myself in 2002 \cite{HP1}.  This led us to study Dirac cohomology of Harish-Chandra modules \cite{HP2} and it became a very useful tool in representation theory.
Kostant  extended Vogan's conjecture to the more general setting of the cubic Dirac operator \cite{Ko2}. 

In the formulation of central problems in Langlands program stable conjugacy plays a pivotal role. The theory of endoscopy investigates
 the difference between orbital integrals over ordinary and stable conjugacy classes.  As Dirac cohomology of a Harish-Chandra
 module determines its K-character \cite{HPZ}, it corresponds to the dual object of the orbital integrals on elliptic elements.  By using
Dirac cohomology,  the transfer factor
is expressed as the difference of characters for the even and odd parts of the spin module, or Dirac index of the  trivial representation \cite{H}. 

The transfer factor of the Adams lifting
is the difference of two irreducible components of the oscillator representation. This transfer factor
is equal to the symplectic Dirac index of the trivial representation for the Lie superalgebra $\fro\frsp(1|2n)$.
This is analogous to the fact
that the endoscopic transfer factor
as the difference of characters for the even and odd parts of the spin module is equal to Dirac index of the  trivial representation \cite{H}. 
The corresponding Vogan's conjecture for symplectic Dirac operators of Lie superalgebras was established by Pandzic and 
myself in 2005
\cite{HP3}.  Let  $Mp(2n,\bbR)$ denote the nonlinear two-fold covering
group of the symplectic group $Sp(2n,\R)$.
By formulating the Adams lifting in the framework of symplectic Dirac cohomology, we obtain
direct lifting of characters  from the linear
group $Sp(2n,\bbR)$ to the nonlinear group $Mp(2n,\bbR)$.
This lifting gives a bijection between the discrete series of $Sp(2n,\bbR)$ and the genuine discrete series of $Mp(2n,\bbR)$.

%%%%%%%%%%%%%%%%%%%%%%%%%%%%%%
\section{Preliminaries on symplectic Dirac cohomology for Lie superalgebras}
%%%%%%%%%%%%%%%%%%%%%%%%%%%%%%%%%%%

Let $\frg=\frg_0\oplus\frg_1$ be a Lie superalgebra.
If $\frg_0$ is reductive, then the adjoint action of
$\frg_0$ on $\frg_1$ is completely reducible. In this case, $\frg$
is called a classical Lie superalgebra.  If a
classical Lie superalgebra $\frg$ is also of Riemannian type, i.e.,
it has nondegenerate supersymmetric invariant bilinear form, then
$\frg$ is called a basic classical Lie superalgebra.  Kac
classified all simple Lie superalgebras \cite{Kac}.  Besides ordinary Lie
algebras the simple basic classical Lie superalgebras list as
 $A(m,n)$, $B(m,n)$, $C(n)$, $D(m,n)$, $D(2,1,\alpha)$,
$F(4)$ and $G(3)$ in \cite{Kac}.

Now we fix a classical complex Lie superalgebra
$$
\frg=\frg_0\oplus\frg_1
$$
with bracket $[\cdot\,,\cdot]$.  We assume that $\frg$ is of Riemannian type
(namely $\frg$ is a basic classical Lie superalgebra), i.e.,
there exists a nondegenerate supersymmetric invariant bilinear 
form $B$ on $\frg$. Supersymmetry means that $B$ is symmetric on $\frg_0$
and skew-symmetric on $\frg_1$, and $\frg_0$ and $\frg_1$ are orthogonal.
Invariance means that
$$
B([X,Y],Z)=B(X,[Y,Z]),
$$
for all $X,Y,Z\in\frg$. 

We fix an orthonormal basis $W_k$ for $\frg_0$ with respect to $B$.
Furthermore, we fix a pair of complementary Lagrangean subspaces of
$\frg_1$, and bases $\del_i$, $x_i$ in them, such that
$$
B(\del_i,x_j)=\frac{1}{2}\delta_{ij}.
$$
This notation is chosen so that $\del_i$ and $x_i$ are generators for
the Weyl algebra $W(\frg_1)$, with only nontrivial commutation relations
being
$$
[\del_i,x_i]_{W(\frg_1)} = 1.
$$
We put the subscript $W(\frg_1)$ to the commutators in $W(\frg_1)$ to distinguish
them from the (totally different) bracket in $\frg$. We see that $W(\frg_1)$
gets identified with the algebra of differential operators with polynomial
coefficients in the $x_i$'s, with $\del_i$ corresponding to
$\frac{\partial}{\partial x_i}$.

Note that if we take
$$
\del_1,\dots,\del_n,x_1,\dots,x_n
$$
for a basis of $\frg_1$, then the dual basis (with respect to $B$) is
$$
2x_1,\dots,2x_n,-2\del_1,\dots,-2\del_n.
$$
Here we say that a basis $f_i$ is dual to the basis $e_i$ if $B(e_i,f_j)=\delta_{ij}$;
namely, for super spaces it is the identification $V\otimes V^*=\Hom(V,V)$ that
involves no signs.

In view of this, the Casimir element of $\frg$ is defined as
$$
\Omega_\frg = \sum_kW_k^2+2\sum_i(x_i\del_i-\del_i x_i).
$$
It is easily checked to be an element of the center $Z(\frg)$ of the
enveloping algebra $U(\frg)$ of $\frg$. Using the relation
$\del_i x_i+x_i\del_i=[\del_i,x_i]$ in $U(\frg)$, one can also write it as
$$
\Omega_\frg = \sum_kW_k^2+4\sum_ix_i\del_i-2\sum_i[\del_i, x_i].
$$
Furthermore, it is independent of the choice of basis: if $e_j$ is any basis
of $\frg$, with dual basis $f_j$ with respect to $B$, then $\Omega_\frg=\sum f_j e_j$.

The action of $\frg_0$ on $\frg_1$ via the bracket defines a map
$$
\nu:\frg_0\lra \frs\frp(\frg_1).
$$
On the other hand, $\frs\frp(\frg_1)$ can be mapped into $W(\frg_1)$
as follows. First we note that the symmetrization map
$\sigma:S(\frg_1)\lra W(\frg_1)$ is a linear isomorphism. Consider the
action of $\sigma(S^2(\frg_1))$ on $\frg_1\subset W(\frg_1)$ by commutators in
$W(\frg_1)$. One readily checks that in our basis $(\del_i,x_j)$, the
commutators with $\sigma(x_ix_j)=x_ix_j$,
$\sigma(\del_i\del_j)=\del_i\del_j$ and $\sigma(\del_ix_j)$ correspond to
the following matrices ($E_{k\;l}$ is the matrix unit, having the $kl$ entry
equal to $1$ and other entries equal to $0$):

\begin{align*}
 \sigma(x_ix_j) & \longleftrightarrow  -E_{n+i\; j}-E_{n+j\; i}; \\
\sigma(\del_i\del_j) & \longleftrightarrow \;\;\; E_{i\; n+j}+E_{j\; n+i}; \\
 \sigma(\del_ix_j) & \longleftrightarrow  -E_{i\; j}+E_{n+j\; n+i}.
\end{align*}

For the last one, note that
$\sigma(\del_i x_j)=\frac{1}{2}(\del_ix_j+x_j\del_i)=\del_ix_j-\frac{1}{2}\delta_{ij}$
acts in the same way as $\del_ix_j$. Since $\frs\frp(\frg_1)$ consists of block matrices
of the form
$$
\begin{pmatrix}
 A & B \cr C & -^tA
\end{pmatrix}
$$
where $A$ is an arbitrary $n\times n$ matrix and $B$ and $C$ are symmetric $n\times n$
matrices, we see that $\sigma(S^2(\frg_1))$ is a Lie subalgebra of $W(\frg_1)$ isomorphic
to $\frs\frp(\frg_1)$ via the isomorphism described above.

Combining with the map $\nu$ defined above, we get a Lie algebra morphism
$$
\alpha:\frg_0 \lra W(\frg_1).
$$
Compared to \cite{Ko}, our $\alpha$ is his $\nu_*$ followed by the symmetrization map.

Now $\alpha$ gives rise to a diagonal embedding
$$
\frg_0\lra U(\frg)\otimes W(\frg_1),
$$
given by
$$
 X\mapsto X\otimes 1+1\otimes\alpha(X).
$$
We denote the image of this map by $\frg_{0\Delta}$; this is a diagonal copy of $\frg_0$.
We denote by $U(\frg_{0\Delta})$, $Z(\frg_{0\Delta})$ the corresponding images of
$U(\frg_0)$ and its center in $U(\frg)\otimes W(\frg_1)$. We will be particularly
interested in the image of the Casimir element $\Omega_{\frg_0}=\sum_k W_k^2$.
It is
$$
\Omega_{\frg_{0\Delta}}=\sum_k (W_k^2\otimes 1+2W_k\otimes\alpha(W_k)+1\otimes\alpha(W_k)^2).
$$
Kostant \cite{Ko} has shown $\alpha(\Omega_{\frg_0})=\sum_k \alpha(W_k)^2$ is a constant
which we denote by $C$
($C$ is equal to $1/8$ of the trace of $\Omega_{\frg_0}$ on $\frg_1$).

Kostant actually showed the following: let $\frg_0$ be a Lie algebra, with a nonsingular
invariant symmetric bilinear form $B$. Let $\frg_1$ be a vector space with
nonsigular alternating bilinear form $B$. Suppose $\frg_1$ is a symplectic representation
of $\frg_0$,
i.e., there is a Lie algebra map $\frg_0\to\frs\frp(\frg_1)$.
Then $\frg =\frg_0\oplus\frg_1$ has a structure of a Lie superalgebra of Riemannian type
compatible with the
given data if and only if $\alpha(\Omega_{\frg_0})$ is a constant. 

Recall that the symplectic Dirac operator for a basic classical Lie superalgebra
corresponding to the decomposition $\frg=\frg_0\oplus\frg_1$ is defined
analogously to the Dirac operator for the Cartan decomposition of a real reductive Lie algebra.
It is an element $D$ of $U(\frg)\otimes W(\frg_1)$ given by
$$
 D=2\sum_i(\partial_i\otimes x_i-x_i\otimes \partial_i).
$$
Then $D$ is easily seen to be independent of the choice of basis, in the sense that if we
 take any basis $e_j$ for $\frg_1$ with dual basis $f_j$ with respect to $B$, then
$D=\sum e_j\otimes f_j$. This also implies that $D$ is $\frg_0$-invariant.
We also have a formula for $D^2$, which is analogous to Parthasarathy's formula for the square of
the Dirac operator attached to the Cartan decomposition for a real reductive Lie algebra:
$$
D^2= -\Omega_\frg \otimes 1 + \Omega_{\frg_{0\Delta}} - C,
$$
where $C$ is the constant described above.

Recall that a Cartan subalgebra of $\frg$ reduces to the Cartan subalgebra of
the even part $\frg_0$.
Let $\frh_0$ be a Cartan subalgebra of $\frg_0$.
We denote by $\Delta$, $\Delta_0$ and $\Delta_1$ the
sets of all roots, even roots and odd roots respectively.
Let $\frb_0$ be a Borel subalgebra of $\frg_0$, containing
$\frh_0$.  We fix a Borel subalgebra  $\frb=\frb_0\oplus\frb_1$
of $\frg$.  Then there exist subalgebras $\frn^+$ and
$\frn^-$ such that
$$\frg=\frn^+\oplus\frh_0\oplus\frn^-, \ \ \frb=\frh_0\oplus\frn^+$$
and
$$ [\frh_0,\frn^+]\subset\frn^+, \ \ [\frh_0,\frn^-]\subset\frn^-.$$
Denote by $\Delta^+$, $\Delta_0^+$ and $\Delta_1^+$ the subsets
of positive roots in the sets $\Delta$, $\Delta_0$ and $\Delta_1$
respectively.  Set
$$\rho_0=\frac{1}{2}\sum_{\alpha\in\Delta_0^+}\alpha,\ \
\rho_1=\frac{1}{2}\sum_{\alpha\in\Delta_1^+}\alpha,
\ \hbox{and}\ \rho=\rho_0-\rho_1. $$

Let $Z(\frg)$ be the center of the enveloping superalgebra $U(\frg)$.
An element $z\in Z(\frg)$ can be uniquely written in the form:
$$z=u_z+\sum_i u_i^+u_i^0u_i^-, \
\hbox{where}\ u_z, u_i^0\in U(\frh_0)\ \hbox{and}\
u_i^\pm\in\frn^\pm U(\frn^\pm). $$
The map $z\mapsto u_z$ gives a monomorphism
$$\beta\colon Z(\frg)\ra U(\frh_0)=S(\frh_0) (=\Bbb C[\frh_0^*]).$$
Let $\tau$ be the automorphism of $\Bbb C[\frh_0^*]$
defined by $(\tau P)(\lambda)=P(\lambda-\rho)$.
The composition $\gamma=\tau\circ\beta\colon Z(\frg)\ra S(\frh_0)$
is called the Harish-Chandra monomorphism.

Let $W$ be the Weyl group of $\frg_0$.  Then $\gamma(Z(\frg))$
is a subalgebra of $S(\frh_0)^W$.  Moreover, the fields of
fractions of $\gamma(Z(\frg))$ and $S(\frh_0)^W$ coincide.

Let $V=V_0\oplus V_1$ be a $\Bbb Z_2$-graded vector space
over $\Bbb C$.  A linear representation $\pi$ of a
Lie superalgebra $\frg=\frg_0\oplus\frg_1$ in $V$ is a
homomorphism from $\frg$ into the superalgebra $\End(V)$
$$\pi\colon \frg\rightarrow \End(V)=\End_0(V)\oplus\End_1(V).$$
For brevity we often say that $V$ is a $\frg$-module.

Let $\lambda\in\frh^*_0$ be a linear functional on $\frh_0$.  We say that a $\frg$-module $V$ has infinitesimal character
$\lambda$ if $Z(\frg)$ acts on $V$
via character 
$$\chi_\lambda\colon Z(\frg)\ra \Bbb C$$ defined
by $\chi_\lambda(z)=\gamma(z)(\lambda)$.  (To define $\gamma(z)(\lambda)$, we 
identify $S(\frh_0)$ with polynomial functions on $\frh_0^*$.)
 Clearly, $\chi_\lambda=\chi_{w\lambda}$ for $w\in W$.
If $V$ is a highest weight $\frg$-module with highest weight
$\Lambda$, then the infinitesimal character of $V$
is $\Lambda+\rho$, or any element in $W\cdot (\Lambda+\rho)$.

Let $M(\frg_1)$ be the
Weil representation for the Weyl algebra $W(\frg_1)$.
It is naturally a $\frsp(\frg_1)$-module (and $\frg_0$-module) 
by the embedding of $\frsp(\frg_1)$ into $W(\frg_1)$ (and the homomorphism
from $\frg_0$ to $\frsp(\frg_1)$).

\begin{defi}
Let $V$ be a $\frg$-module.  The Dirac cohomology
$H_D(V)$ of $V$ is defined to be the $\frg_0$-module
$$
\ker D/\ker D\cap \im D.
$$
\end{defi}

\begin{thm} \cite{HP3} Let $\frg$ be a basic classical Lie superalgebra and
$V$ a $\frg$-module with infinitesimal character $\chi$.
If the Dirac cohomology $H_D(V)$
contains a nonzero $\frg_0$-module with infinitesimal character
$\lambda\in \frh_0^*$, then $\chi$ is determined by the $W$-orbit
of $\lambda$.
More precisely, there exists a homomorphism of
$Z(\frg)$ into $Z(\frg_0)$ given by the following
commuting diagram such that $\chi(z)=\chi_\lambda(\zeta(z))$
for $z\in Z(\frg)$:
$$
\CD
Z(\frg)@>\zeta >>Z(\frg_0) \\
@V{\hbox{H.C. hom}}VV                 @VV\hbox{H.C. isom}V     \\
S(\frh_0)^W@>{\hbox{id}}>>S(\frh_0)^W
\endCD
$$
where the bottom horizontal map is identity and the
two vertical maps are the Harish-Chandra monomorphism and
isomorphism respectively.

\end{thm}

If we identify $W(\frg_1)$ as the algebra of differential
operators with polynomial coefficients in the $x_i$'s ($i=1,\cdots, n$), then
the Weil representation $M(\frg_1)$ can be identified as the polynomial algebra
$\Bbb C[x_1,\cdots ,x_n]$.  We write $M^+(\frg_1)$
and $M^-(\frg_1)$ for the $\frsp(\frg_1)$-submodules of $M(\frg_1)$ spanned by homogeneous
polynomials of even and odd degrees respectively.
Clearly,
$$D\colon V\otimes M^\pm(\frg_1)\ra V\otimes M^\mp(\frg_1).$$
Now we consider
\[
D^+:V\otimes M^+\rightarrow V\otimes M^- \text{ and } D^-:V\otimes M^-\rightarrow V\otimes M^+
\]
defined by the restrictions of $D$.
It follows that
$$H_D(V)=H_D^+(V)\oplus H_D^-(V)$$
 and
\[
V\otimes M^+-V\otimes M^-=H_D^+(V)-H_D^-(V).
\]
In terms of $\frg_0$-characters, this reads
$$
\ch V(\ch M^+-\ch M^-)=\ch H_D^+(V)-\ch H_D^-(V).
$$
In particular, set $V=1\!\!1$ (the trivial representation), we have
$$\ch M^+-\ch M^-=\ch H_D^+(1\!\!1)-\ch H_D^-(1\!\!1).$$ 
This gives a $\frg_0$-character formula for $V$ analougous to the Weyl character formula.

\begin{thm} Let $\frg$ be a basic classical Lie superalgebra and
$V$ a finite-dimensional $\frg$-module.  Then we have the following formula for $\frg_0$-characters 
$$
\ch V=\frac{\ch H_D^+(V)-\ch H_D^-(V)}{\ch M^+-\ch M^-}=\frac{\ch H_D^+(V)-\ch H_D^-(V)}{\ch H_D^+(1\!\!1)-\ch H_D^-(1\!\!1)}.
$$
\end{thm}

%%%%%%%%%%%%%%%%%%%%%%%%%%%%%%%%%%
\section{Rittenberg-Scheunert correspondence}
%%%%%%%%%%%%%%%%%%%%%%%%%%%%%%%%%%%

We now discuss Rittenberg-Scheunert correspondence for representations of $\fro\frsp(1|2n)$ and  $\fro(2n+1)$.
 The roots of the complex basic Lie algebra $\frg=\fro\frsp(1|2n)$ fall into the two
 subsets of even roots
$$\Delta_0=\{\pm e_i\pm e_j;  1\leq i< j\leq n\}\cup \{\pm 2e_i; 1\leq i \leq n\}, $$
and odd roots
$$ \Delta_1=\{\pm e_i; 1\leq i \leq n\}.$$
We select simple roots
$$e_1-e_2, e_2-e_3,\cdots, e_{n-1}-e_n, e_n.$$
Then the positive  even roots are
$$ e_i\pm e_j \text{ with } 1\leq i< j\leq n, \  2e_i \text{ with }1\leq i \leq n, $$
and positive odd roots are
$$  e_i \text{ with } 1\leq i \leq n.$$
The corresponding fundamental weights are
$$\omega_i=e_1+\cdots +e_i  \ (i=1,\ldots, n-1) \text{ and } \omega_n={1\over 2}(e_1+\cdots +e_n).$$
Note that the even roots $\Delta_0$ are the roots of $\frsp(2n)$ and 
$\omega_1, \ldots, \omega_{n-1},2\omega_n$ are the fundamental weights for $\frsp(2n)$.

Let $\rho_0$ and $\rho_1$ denote the half sum of the positive even roots, respectively odd roots.
Then
$$\rho_0=\omega_1+ \cdots +\omega_{n-1}+ 2\omega_n$$
and 
$$\rho_1=\omega_n.$$
We set
$\rho=\rho_0-\rho_1$.
Then 
$$\rho=\omega_1+ \cdots +\omega_{n-1}+ \omega_n.$$
Writing in  coordinates, 
$$\rho_0=(n,\ldots,1),\ \rho_1=({1\over 2},\ldots,{1\over 2}), \ \rho=(n-{1\over 2}, n-{3\over 2},\ldots, {1\over 2}).$$
The irreducible finite-dimensional representations of 
$\frg=\fro\frsp(1|2n)$ are uniquely determined by their highest weights 
of the form
$$\Lambda=p_1\omega_1+\cdots +p_n\omega_n,$$
with integers $p_i\geq 0$ and $p_n$ even.  These are also
the highest weights of finite-dimensional representations of $\frg_0=\frsp(2n)$.

The setting for the complex simple Lie algebra $\fro(2n+1)$ is similar.
The root system of $\fro(2n+1)$ consists of long roots
$$\pm e_i\pm e_j \text{ with  } 1\leq i< j\leq n, $$
and short roots
$$ \pm e_i \text{ with  } 1\leq i \leq n.$$
We select simple roots
$$e_1-e_2, e_2-e_3,\cdots, e_{n-1}-e_n, e_n.$$
Then the positive long roots are
$$ e_i\pm e_j \text{ with  }  1\leq i< j\leq n, $$
and positive short roots are
$$  e_i \text{ with  }  1\leq i \leq n.$$
The corresponding fundamental weights are
$$\omega_i=e_1+\cdots +e_i \  (i=1,\ldots, n-1) \text{ and } \omega_n={1\over 2}(e_1+\cdots +e_n).$$

Let $\rho$ denote the sum of the positive  roots.
Then
$$\rho=\omega_1+ \cdots +\omega_{n-1}+ \omega_n$$

The irreducible finite-dimensional representations of 
$\frg=\fro(2n+1)$ are uniquely determined by their highest weights 
of the form
$$\Lambda=p_1\omega_1+\cdots +p_n\omega_n,$$
with integers $p_i\geq 0$ ($i=1,\ldots, n$). For non-spinorial representations,  $p_n$ must be even,
which are exactly the same condition as for 
the highest weights of finite-dimensional representations of $\frg_0=\frsp(2n)$ or  $\frg=\fro\frsp(1|2n)$.

\begin{thm}\cite[Sect. 2]{RS}
The highest weights of both graded irreducible representations
of $\fro\frsp(1|2n)$ and the non-spinorial irreducible representations of $\fro(2n+1)$
are the same, taking the form
$$\Lambda=p_1\omega_1+\cdots +p_n\omega_n,$$
with integers $p_i\geq 0$ and $p_n$ even.  
If a  graded irreducible representation
of $\fro\frsp(1|2n)$ and a non-spinorial irreducible representation of $\fro(2n+1)$
have the same highest weight, then the multiplicity of any weight is the same for both representations.
\end{thm}

\begin{thm}\cite[Note added in proof]{RS}
If a  graded irreducible representation
of $\fro\frsp(1|2n)$ and a non-spinorial irreducible representation of $\fro(2n+1)$
have the same highest weight 
$$\Lambda=p_1\omega_1+\cdots +p_n\omega_n,$$
then they have the same corresponding parameters  $\Lambda+\rho$ as their infinitesimal characters.
\end{thm}

The Rittenberg-Scheunert correspondence for $\fro\frsp(1|2n)$ and $\fro(2n+1)$ was extended to quantum groups
$U_q(\fro\frsp(1|2n))$ and $U_{-q}(\fro(2n+1))$ by Ruibin Zhang \cite{Z}.

%%%%%%%%%%%%%%%%%%%%%%%%%%%%%%%%%%
\section{Lifting of characters to metaplectic groups}
%%%%%%%%%%%%%%%%%%%%%%%%%%%%%%%%%%%

Adams defined a lifting of characters between special orthogonal groups $SO(2n+1)$ over $\R$
and the nonlinear metaplectic groups $Mp(2n,\R)$ \cite{A}.  It is closely related to both 
endoscopy and theta-lifting.   If $\pi$ is an irreducible representation of $SO(p,q)$ with $p+q=2n+1$, 
then $\pi$ has a non-zero theta-lift $\pi'$ to a genuine representation of $Mp(2n,\R)$.
The Adams lifting amounts to addressing the relation between
the characters of $\pi$ and $\pi'$.
The main idea of \cite{A} is interpreting $G=SO(n+1,n)$ as an endoscopic group for $G'=Mp(2n,\R)$ 
and the corresponding character lifting is best illustrated by
the examples of the discrete series.  Let $\pi_\lambda$ be a discrete series representation of $G=SO(n+1,n)$
with Harish-Chandra parameter
$$\lambda=(a_1,\ldots, a_k, b_1,\ldots b_l) $$
with $a_i,b_j\in \Z+\frac{1}{2}, a_1>\cdots>a_k>0, b_1>\cdots>b_l>0$. 
Under the theta-lifting, $\pi_\lambda$ corresponds to a discrete series $\pi_{\lambda'}$ of $Mp(2n,\R)$ 
with Harish-Chandra parameter
$$\lambda'=(a_1,\ldots, a_k,-b_l,\ldots, -b_1).$$
Denote by $W$ the Weyl group of $G$ or $G'$ (type $B_n$ or $C_n$),  $W_K$ the Weyl group 
of $K=S(O(n+1)\times O(n))$ and  $W_{K'}$ the Weyl group 
of $K'=U(n)$.
We consider the two stable sums of discrete series representations
$${\bar \pi_\lambda}=\sum_{w\in W_K\backslash W}\pi_{w\lambda}$$
and
$${\bar \pi_{\lambda'}}=\sum_{w\in W_{K'}\backslash W}\pi_{w\lambda'}.$$ 
Fix an isomorphism of a compact Cartan subgroup $T$ of $G=SO(n+1,n)$ 
and a compact Cartan subgroup $T'$ of $Sp(2n,\R)$.  Let $\widetilde{T'}$ be the covering group of $T'$ in $G'=Mp(2n,\R)$. 
The character of ${\bar \pi_\lambda}$ on $T$ is 
$$\Theta_{SO}(\lambda)(t)=\frac{\sum_{w\in  W}\sgn(w)e^{w\lambda}(t)} {D_{SO}(t)} \ (t\in T), $$
and similarly, the character of ${\bar \pi_{\lambda'}}$ on $\widetilde{T'}$ is
$$\Theta_{Sp}(\lambda')(t')=\frac{\sum_{w\in  W}\sgn(w)e^{w\lambda'}(t')} {D_{Sp}(t')} \ (t'\in \widetilde{T'}).$$
Then the two characters are related by multiplying by the transfer factor 
$$\Theta_{Sp}(\lambda')(t')=\Theta_{SO}(\lambda)(t)\Phi(t'),$$
where the transfer factor $\Phi(t') =\frac{D_{SO}(t)}{D_{Sp}(t')}$
 is precisely
the difference of two irreducible components of the oscillator representation (up to a sign)
$$\pm \Phi(t')=\Omega_{even}(t')-\Omega_{odd}(t').$$
Note that  the stable character $\Theta_{SO}(\lambda)(t)$ on $T$ coincides with the 
character of the finite-dimensional representation with infinitesimal character $\lambda+\rho$,
which equals to the character of the finite-dimensional representation of $\fro\frsp(1|2n)$ with the same
infinitesimal character under 
the Rittenberg-Scheunert correspondence described in the previous section.
The transfer factor is equal to the symplectic Driac index of the trivial repesentation
$$\Phi=\Omega_{even}-\Omega_{odd}=\ch H_D^+(1\!\!1)-\ch H_D^-(1\!\!1),$$
and the virture character
$$\Theta_{Sp}(\lambda')(t')=\Theta_{SO}(\lambda)(t)\Phi(t')$$
is the symplectic Dirac index of the corresponding finite-dimensional representation of $\fro\frsp(1|2n)$.

In light of the above discussion, we may define a direct lifting of characters of $Sp(2n,\R)$ to $Mp(2n,\R)$.
In a sense, we interpret $Sp(2n,\R)$ as an endoscopic group for $Mp(2n,\R)$.
We retain the notation of the previous section.  Recall that $\omega_n=\rho_1=({1/2},\ldots, {1/2})$.
An invariant eigendistribution $\Theta$
on $Sp(2n,\R)$ can be expressed on a Cartan subgroup by
$$\Theta=\frac{\sum_{w\in W}a_we^{w\lambda}}{D_{SP}}, \text{ with }a_w\in \bbC.$$ 
  We may assume $\lambda$ to be dominant by reorganizing the coefficients $a_w$.
Set 
$$\Theta'=\frac{\sum_{w\in W}a_we^{w(\lambda-\omega_n)}}{D_{SP}}.$$
Then $\Theta'$ is a genuine invariant eigendistribution on $Mp(2n,\R)$. 
We define the lifting $\Gamma$ by
$$\Gamma \colon \Theta \rightarrow \Theta'.$$

If $\lambda=(a_1,\ldots, a_n)$
with $a_i\in \Z, a_1>\cdots>a_n>0$ is the Harish-Chandra parameter of a holomorphic discrete series representation
 for $Sp(2n,\R)$,  then 
$$\lambda'=(a_1,\ldots, a_n)-({1/2},\ldots, { 1/2})=\lambda-\omega_n$$
is  the Harish-Chandra parameter of a holomorphic discrete series representation for $Mp(2n,\R)$.
In general, if $\lambda=(a_1,\ldots, a_k,-b_l,\ldots, -b_1)$ with $a_i,b_j\in \Z, a_1>\cdots>a_k>0, b_1>\cdots>b_l>0$
 is the Harish-Chandra parameter of a discrete series representation for $Sp(2n,\R)$, 
 then
 $$\lambda'=(a_1-{1\over 2},\ldots, a_k-{1\over 2},-b_l+{1\over 2},\ldots, -b_1+{1\over 2})$$
is the Harish-Chandra parameter for a discrete series representation for $Mp(2n,\R)$.
In other words, if $w\in W$ makes $w\lambda$ dominant, then we set the corresponding
$\lambda'$ to be
$$\lambda'=\lambda-w^{-1}\omega_n.$$
Let $\pi_\lambda$ be a discrete series representation of $Sp(2n,\R)$ with the
Harish-Chandra parameter $\lambda$.  Let $\pi_{\lambda'}$ be the corresponding
genuine discrete series representation of $Mp(2n,\R)$ with the
Harish-Chandra parameter $\lambda'$.  Then $\Gamma(\Theta_{\pi_\lambda})=\Theta_{\pi_{\lambda'}}$.
Consequently, we obtain the following theorem.

\begin{thm}
The map $\Gamma$ is a bijection between stably invariant eigendistributions on $Sp(2n,\R)$ and
genuine stably invariant eigendistributions on $Mp(2n,\R)$. It restricts to a bijection between discrete
series of $Sp(2n,\R)$ and genuine discrete
series of $Mp(2n,\R)$.
\end{thm}

%%%%%%%%%%%%%%%%%%%%%%%%%%%%

\end{document}